\documentclass{article}

\usepackage{arxiv}
\usepackage{amssymb}
\usepackage{amsmath}
\usepackage{empheq}
\usepackage{amscd}
\usepackage{graphicx}
\usepackage{graphics}
\usepackage[noadjust]{cite}
\usepackage{amsthm}
\usepackage{caption}
\usepackage{multirow}
\usepackage{array}
\usepackage{rotating}

\usepackage{indentfirst}
\usepackage{tikz}
\usepackage{calc}
\usepackage{epsfig}
\usepackage{natbib}
\usepackage[makeroom]{cancel}
\usepackage{multicol}
\usepackage{csquotes}
\usepackage{algorithm}
\usepackage{algorithmic}
\usepackage{enumitem}
\usepackage{hyperref}
\usepackage{url}
\usepackage{bbm}  
\hypersetup{colorlinks,linkcolor={blue},citecolor={blue},urlcolor={blue}}

\usepackage{timet}
\usepackage{graphicx}
\usepackage{multirow}
\usepackage{multicol}
\usepackage{lipsum}
\usepackage{timet}
\usepackage{epsfig}
\usepackage{amsmath}
\usepackage{amsfonts}
\usepackage{amssymb}
\usepackage{color}
\numberwithin{equation}{section}
\usepackage{caption}
\usepackage{float}
\usepackage{subcaption}
\usepackage{graphics}
\usepackage{xcolor,graphicx}
\usepackage{csquotes}
\usepackage{mathtools}
\usepackage{optidef}

\DeclareMathOperator*{\minimize}{minimize}
\newcommand{\vertiii}[1]{{\left\vert\kern-0.25ex\left\vert\kern-0.25ex\left\vert #1 
    \right\vert\kern-0.25ex\right\vert\kern-0.25ex\right\vert}}

\begin{document}
\title{Localized Wavefield Inversion (LWI): an Adaptation of Multi-Block ADMM for Localized FWI}

\author{\href{http://orcid.org/0000-0003-1805-1132}{\includegraphics[scale=0.06]{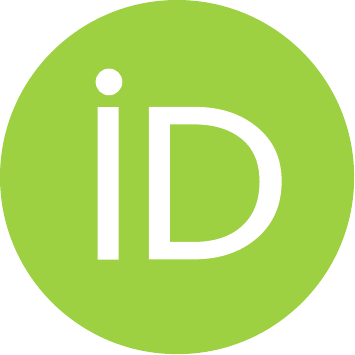}\hspace{1mm}Hossein S. Aghamiry} \\
  University Cote d'Azur - CNRS - IRD - OCA, Geoazur, Valbonne, France. 
  \texttt{aghamiry@geoazur.unice.fr}
\And
  \href{https://orcid.org/0000-0002-9879-2944}{\includegraphics[scale=0.06]{orcid.pdf}\hspace{1mm}Ali Gholami} \\
  Institute of Geophysics, University of Tehran, Tehran, Iran.
  \texttt{agholami@ut.ac.ir} \\ 
  \And
\href{http://orcid.org/0000-0002-4981-4967}{\includegraphics[scale=0.06]{orcid.pdf}\hspace{1mm}St\'ephane Operto} \\ 
  University Cote d'Azur - CNRS - IRD - OCA, Geoazur, Valbonne, France. 
  \texttt{operto@geoazur.unice.fr}
\And
\href{http://orcid.org/0000-0003-2775-7647}{\includegraphics[scale=0.06]{orcid.pdf}\hspace{1mm}Alison Malcolm} \\ 
  Department of Earth Sciences, Memorial University of Newfoundland, St John’s NL A1B 3X5, Canada. 
  \texttt{amalcolm@mun.ca} 
  }

\renewcommand{\shorttitle}{Localized Wavefield Inversion (LWI), Aghamiry et al.}

\maketitle

\begin{abstract}
Full-waveform inversion (FWI) is a high-resolution and computationally intensive imaging technique to reconstruct unknown parameters in the computational model in which the waves propagate; however, an accurate model of only part of this medium is required for some applications. To decrease the computational burden of such problems, target-oriented FWI was proposed where the redatumed data on the part of the medium or localized solvers for the wave equation are used.\\
On the other hand, the classical formulation of FWI suffers from non-linearity and ill-posedness, which makes FWI sensitive to the initial model, the low-frequency content of the data, and limited illumination. In this study, we propose a localized version of the alternating direction method of multipliers (ADMM)-based FWI method, which was proposed to solve these problems in classical FWI. In our localized FWI or LWI, the medium is decomposed into a few subdomains, where some of them are updated, and the others are kept fixed based on an adaptation of multi-block ADMM, which is a powerful algorithm for solving inverse problems with decomposition and block separability. Numerical tests on the Marmousi model for a time-lapse application confirm the computational efficiency and robustness against background velocity model errors. 
\end{abstract}

\section{INTRODUCTION}
\begin{figure}[t!]
\center
\includegraphics[width=0.4\columnwidth,clip=true,trim=9.2cm 7.8cm 12.5cm 5.2cm]{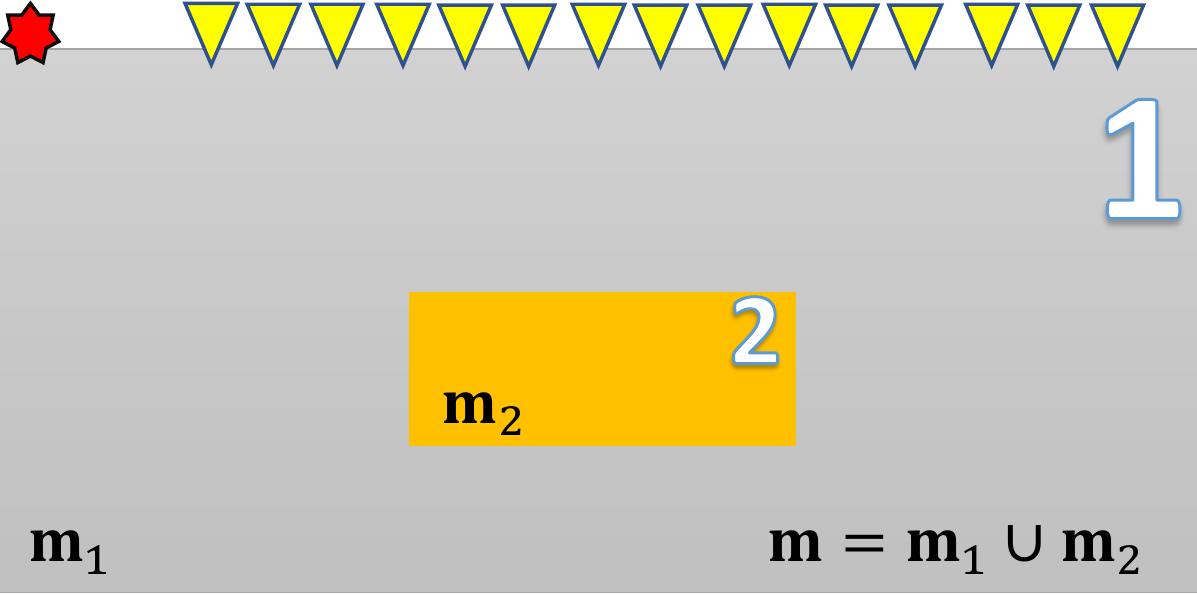}
\caption{The schematic of localized FWI. The goal is to estimate the model parameters in sub-domain $\mathbbm{2}$ when an approximation of sub-domain $\mathbbm{1}$ ($\bold{m}_1$) is available. The red star and yellow rectangles show the source and receivers, respectively. Note that, the union of $\bold{m}_1$ and $\bold{m}_2$ is $\bold{m}$.}
\label{fig:FIG1}
\end{figure}

Seismic full-waveform inversion (FWI) is an accurate imaging technique that can reconstruct the high-resolution elastic properties of the medium \citep{Pratt_1998_GNF}.
On the other hand, it suffers from non-linearity, ill-posedness, and a huge computational burden, all of which need to be addressed to make FWI a practical tool for imaging in different  applications \citep{Virieux_2009_OFW}.\\
FWI generally estimates model parameters in the whole computational domain in which the seismic waves propagate, which is computationally expensive and not always necessary. 
For some specific target-oriented applications, often related to time-lapse imaging like the monitoring of reservoirs during hydrocarbon production or CO$_2$ injection \citep{Willemsen_2017_PWA}, or other applications like salt boundary inversion \citep{Lewis_2012_ALS} or geothermal exploitation, one only seeks to update model parameters or their evolution over time within a localized region of interest. Such context is also encountered in earthquake seismology for lithospheric imaging from teleseismic events \citep{Monteiller_2013_NOA,Beller_2018_LAS,Wang_2022_FWI}, or at the global scale for deep Earth investigations \citep{Masson_2017_BTL}.\\
Target-oriented FWI \citep{Malcolm_2016_R4F,Willemsen_2016_ANE,Yuan_2017_E3L} is a computationally-efficient method for such kinds of problems. Fig. \ref{fig:FIG1} illustrates the concept: instead of updating the model parameters in the whole domain ($\bold{m}$), we seek to update them only in a small subdomain, indicated by $\mathbbm{2}$ in Fig. \ref{fig:FIG1}, when a good approximation of the surrounding medium, indicated by $\mathbbm{1}$, is available.\\ 
The conventional methods for target-oriented FWI are based on data redatuming \citep{Yang_2012_TTW,Vanderneut_2015_OGF} or local wave-equation solvers \citep{Robertson_2000_EMC,Van_2007_EWF,Borisov_2015_EME,Vanderneut_2015_OGF,Willemsen_2016_ANE}. 
In the first method, the sources and receivers are moved to $\mathbbm{2}$, hence, reducing the size of the problem and decreasing the computation time. 
The redatuming needs an integral over the physical source and receivers  \citep{Wiggins_1984_KIE}; however, the finite sampling of the integrals,  limited acquisition coverage, and  inaccurate velocity models introduce artifacts in the virtual redatumed data \citep{Mulder_2005_RE,Haffinger_2013_SBF}.\\   
The second category of methods do not modify the data acquisition setup, but they try to reconstruct an accurate local wavefield in $\mathbbm{2}$  \citep{Robertson_2000_EMC,Teng_2003_EBC,Van_2007_EWF,Gillman_2015_ASA,Masson_2017_BTL,Broggini_2017_IBC}. 
\citet{Robertson_2000_EMC} propose the finite-difference injection method where initially the wavefields on the boundary of $\mathbbm{2}$ are stored after a wavefield simulation in the entire medium. After alteration of $\mathbbm{2}$, the wavefields are updated in $\mathbbm{2}$ by using the stored wavefield, whereas the surrounding model is kept fixed. 
This method has a straightforward implementation, and it doesn't require any boundary conditions for $\mathbbm{2}$. 
But, on the other hand, this method doesn't account for the scattering of the updated wavefield occuring in $\mathbbm{1}$ and re-entering $\mathbbm{2}$, which is called high-order long-range interaction by \citet{Robertson_2000_EMC}. 
Inspired by \citet{Van_2007_EWF,Vasmel_2013_IEI}, \citet{Willemsen_2017_PWA} proposed another localized solver method in the frequency domain by applying immersive boundary condition (IBC) to $\mathbbm{2}$. Such IBCs take into account all of the wavefield interactions between altered subdomain $\mathbbm{2}$ and unchanged subdomain $\mathbbm{1}$, but, for each frequency, it requires as many wave-equation solutions as grid points discretizing the edges of the domain of interest for simulating the Green's functions in the entire medium. To keep the computational burden reasonable, \citet{Kumar_2019_ENE} proposed a rank minimization-based framework to compute a low-rank factorized form of the Green's functions, which requires fewer wave-equation solutions.\\
FWI is classically solved based on a reduced approach formulation with local optimization algorithms \citep{Pratt_1998_GNF}, the only affordable approach to the problem considering the size of the data and model spaces, particularly in 3D. Therefore, FWI requires an initial subsurface model that is iteratively improved until convergence. The inversion is rapidly trapped in local minima when the initial model does not predict recorded traveltimes to less than half a period \citep{Virieux_2009_OFW}. Overcoming this issue has attracted some attention over the last decade \citep[ and the references therein]{Gholami_2022_EFW}.\\
Among such methods, one category consists in extending the linear regime of FWI by enlarging the search space in the state-variable domain \citep{Abubakar_2009_FDC,VanLeeuwen_2013_MLM,Aghamiry_2019_IWR} among others. Considering the bi-convexity of the FWI problem, \citet{Aghamiry_2019_IWR} propose using the alternating direction method of multipliers (ADMM) \citep{Boyd_2011_DOS}, which is an efficient method for solving bi-convex optimization problems.
This method, which is called iteratively refined wavefield reconstruction inversion (IR-WRI), is an attempt to solve the FWI problem in full-space, i.e., when the wavefields, model parameters, and Lagrange multipliers are updated independently \citep{Rieder_2021_AAF}, namely all-at-once methods.
IR-WRI updates these three classes of variables in an alternating mode. The method relaxes the requirement to satisfy the wave equation at each iteration to improve the data fit. This amounts to reconstructing wavefields, namely data-assimilated (DA) wavefields, that best jointly fit the observations and satisfy the wave equation in the least-squares sense. Then, the subsurface parameters are updated by minimizing the wave-equation errors generated by the relaxation. Finally, the Lagrange multipliers (dual variables) are updated with a classical dual ascent approach, where they accumulate the remaining errors of the constraints. \\
In this study, we want to develop a localized version of IR-WRI (LWI) to update $\mathbbm{2}$ when a good approximation of $\mathbbm{1}$ is available. 
We first do a block decomposition on FWI to separate $\mathbbm{1}$ and $\mathbbm{2}$. 
By doing so, FWI becomes a constrained optimization problem with four unknowns, wavefields, and model parameters in $\mathbbm{1}$ and $\mathbbm{2}$. This problem can be solved with multi-block ADMM \citep{Liu_2015_MBA}, which updates all four variables in an alternating mode. This method will have a slower convergence rate, but it gets some computational gain from solving smaller systems compared to IR-WRI.   
In LWI, we propose an adaptation of multi-block ADMM, and we only update the model and wavefields in $\mathbbm{2}$ at each iteration, while we update the model and wavefields in $\mathbbm{1}$ only once. The proposed method is straightforward, and it is computationally cheap, but on the other hand, it estimates an approximate wavefield in the updated subdomain $\mathbbm{2}$.   
Our numerical tests show that this effect is negligible in our formulation because we involve the data when we extract the wavefield in $\mathbbm{1}$. We assess the performance of the method using a time-lapse example with the Marmousi model. The results show that the proposed method is robust to errors in the velocity model in $\mathbbm{1}$. In addition, it can update the velocity model in $\mathbbm{2}$ in a reasonable time with a quality close to the case when we apply conventional IR-WRI on the entire domain.  

%
%
%
\section{Theory}
Frequency-domain FWI can be formulated as  
\begin{equation} 
\minimize_{\bold{U},\bold{m}\in \mathcal{M}}~\mathcal{R}(\bold{m}) ~~\text{subject to}
~~\bold{A}(\bold{m})\bold{U}=\bold{B}~~ \text{and}~~ \bold{PU}=\bold{D},\label{main}
\end{equation} 

where $\bold{m} \in \mathbb{C}^{n\times 1}$ is the model parameter vector (squared slowness); $n$ is the number of discretized points of the model; $\bold{A}(\bold{m})\in \mathbb{C}^{n\times n}$ is the Helmholtz operator; $\bold{U} \in \mathbb{C}^{n\times n_s}$, $\bold{B} \in \mathbb{C}^{n\times n_s}$ and $\bold{D} \in \mathbb{C}^{n_r\times n_s}$ denote the wavefield, source term and the recorded data for $n_s$ sources and $n_r$ receivers, respectively; $\bold{P} \in \mathbb{R}^{n_r\times n}$ is the observation operator; $\mathcal{R}(\bold{m})$ is an appropriate regularization function on the model domain; and $\mathcal{M}$ is a convex set defined according to our prior knowledge of $\bold{m}$ \citep{Aghamiry_2019_CRO}. Finally, $\bold{A}(\bold{m})=\Delta+\omega^2 \text{Diag}(\bold{m})$, $\omega$ is the angular frequency, $\Delta$ is the Laplacian, Diag($\bold{x}$) denotes a diagonal matrix with the entries of the vector $\bold{x}$ on its diagonal.
 
Suppose subvectors $\bold{m}_1\in \mathbb{C}^{n_1\times 1}$ and $\bold{m}_2\in \mathbb{C}^{n_2\times 1}$ are the partitions of vector $\bold{m}$ related to $\mathbbm{1}$ and $\mathbbm{2}$, respectively, the wave-equation can be decomposed as \citep{Masson_2017_BTL}
\begin{equation}
\bold{A}(\bold{m})\bold{U}=\bold{A}_1(\bold{m}_1)\bold{U}_1+\bold{A}_2(\bold{m}_2)\bold{U}_2=\bold{B},
\end{equation} 

where $\bold{A}_1(\bold{m}_1) \in \mathbb{C}^{n\times n_1}$ and $\bold{A}_2(\bold{m}_2) \in \mathbb{C}^{n\times n_2}$ are over-determined matrices containing the columns of $\bold{A}(\bold{m})$ related to $\mathbbm{1}$ and $\mathbbm{2}$, respectively, and $n_1+n_2=n$. The same decomposition can be applied on $\bold{P}$ to extract the observation operator in $\mathbbm{1}$ and $\mathbbm{2}$ ($\bold{P}_1$, and $\bold{P}_2$, respectively). 
Accordingly, FWI can be written as 
\begin{equation}
\minimize_{\bold{U}_1,\bold{U}_2,\bold{m}_1, \bold{m}_2\in \mathcal{M},}~~~\mathcal{R}_1(\bold{m}_1)+\mathcal{R}_2(\bold{m}_2)~~~~~~~~
\text{subject to}~~~~
 \begin{cases} \bold{A}_1(\bold{m}_1)\bold{U}_1+\bold{A}_2(\bold{m}_2)\bold{U}_2=\bold{B},\\ \bold{P}_1 \bold{U}_1+\bold{P}_2 \bold{U}_2=\bold{D}.\end{cases}~~~~~~~~~~~~
\label{main_2}
\end{equation}

In this study, we assume that the receivers are only in $\mathbbm{1}$, .i.e. $\bold{P}_2=\emptyset$, which is the conventional setup for localized inversions in exploration seismology. But, generally, there isn't any restriction for the positions of receivers. For example, like teleseismic applications, they can be in $\mathbbm{2}$, .i.e. $\bold{P}_1=\emptyset$. \\
LWI is based on multi-block ADMM \citep{Boyd_2011_DOS,Liu_2015_MBA} for solving \eqref{main_2}. (Multi-block) ADMM finds the saddle point of a combination of classical Lagrangian and a penalty term, namely the augmented Lagrangian (AL), in an alternating mode for primal variables and dual variables as 
\begin{subequations} 
\label{ADMM}
 \begin{align}
\bold{U}_1^{k+1}&= \underset{\bold{U}_1}{\arg\min} ~ \Psi(\bold{U}_1,\bold{U}^k_2,\bold{m}_1^{k},\bold{m}_2^{k},\bold{\hat{B}}^k,\bold{\hat{D}}^k), \label{u1}\\ 
\bold{U}_2^{k+1}&= \underset{\bold{U}_2}{\arg\min} ~ \Psi(\bold{U}^{k+1}_1,\bold{U}_2,\bold{m}_1^{k},\bold{m}_2^{k},\bold{\hat{B}}^k,\bold{\hat{D}}^k),\label{u2}\\ 
\bold{m}_1^{k+1}&= \underset{\bold{m}_1\in \mathcal{M}}{\arg\min} ~ \Psi(\bold{U}^{k+1}_1,\bold{U}^{k+1}_2,\bold{m}_1,\bold{m}_2^{k},\bold{\hat{B}}^k,\bold{\hat{D}}^k), \label{m1}\\
\bold{m}_2^{k+1}&= \underset{\bold{m}_2\in \mathcal{M}}{\arg\min} ~ \Psi(\bold{U}^{k+1}_1,\bold{U}^{k+1}_2,\bold{m}_1^{k+1},\bold{m}_2,\bold{\hat{B}}^k,\bold{\hat{D}}^k),\label{m2}\\
\bold{\hat{B}}^{k+1} &= \bold{\hat{B}}^k  +\bold{B} - \bold{A}_1(\bold{m}_1^{k+1})\bold{U}_1^{k+1}-\bold{A}_2(\bold{m}_2^{k+1})\bold{U}_2^{k+1}, \label{Bk}\\
\bold{\hat{D}}^{k+1} &= \bold{\hat{D}}^k  +\bold{D}- \bold{P}\bold{U}_1^{k+1}\label{Dk}, 
 \end{align}
\end{subequations} 
where
\begin{align}
&\Psi(\bold{U}_1,\bold{U}_2,\bold{m}_1,\bold{m}_2,\bold{\hat{B}}^k,\bold{\hat{D}}^k)=
\mathcal{R}_1(\bold{m}_1)+\mathcal{R}_2(\bold{m}_2)+ 
\|\bold{PU}_1-\bold{D}-\bold{\hat{D}}^k\|_F^2 +\lambda \|\bold{A}_1(\bold{m}_1)\bold{U}_1+\bold{A}_2(\bold{m}_2)\bold{U}_2-\bold{B}-\bold{\hat{B}}^k\|_F^2,\nonumber \label{eqpsi}
\end{align}
is the scaled form of AL, $F$ is the Frobenius norm, $\bullet^k$ is the value of $\bullet$ at iteration $k$, $\lambda>0$ is the penalty parameter assigned to the wave equation constraint, and $\bold{\hat{B}}^{k}\in \mathbb{C}^{n\times n_s}$, $\bold{\hat{D}}^{k}\in \mathbb{C}^{n_r\times n_s}$ are the scaled dual variables. \\
Solving Eqs. \eqref{ADMM} will converge to the solution of \eqref{main}, but with a slower convergence rate compared to the conventional ADMM (IR-WRI). In LWI, we only want to update the model parameters in $\mathbbm{2}$. 
To this end, we modify Eqs. \eqref{ADMM} as follows:\\ 

\textbf{1) Estimating the wavefield in $\mathbbm{1}$:} 
The wavefield in $\mathbbm{1}$ is updated only once, at the first iteration of the inversion of each frequency, because we need it to update the wavefield and model parameters in $\mathbbm{2}$ (Eqs. \eqref{u2} and \eqref{m2}). On the other hand, to extract the wavefield in $\mathbbm{1}$ (Eqs. \eqref{u1}), we need the wavefield and model parameters in $\mathbbm{2}$. At this stage, we can jointly estimate the wavefield in $\mathbbm{1}$ and $\mathbbm{2}$, which is nothing more than the DA wavefield of IR-WRI in the entire of the medium as
\begin{equation} \label{DA0}
\bold{U}^0=[\lambda\bold{A(m}^0)^T\bold{A(m}^0)+\bold{P}^T\bold{P}]^{-1}[\lambda\bold{A(m}^0)^T\bold{B}+\bold{P}^T\bold{D}].
\end{equation} 
Then $\bold{U}^0_1$ and $\bold{U}^0_2$ are easily extracted from $\bold{U}^0$. The Lagrange multipliers ($\bold{\hat{B}}$, and $\bold{\hat{D}}$) are initialized with 0, thus these terms are absent in Eq. \eqref{DA0}.\\
Having computed $\bold{U}^0$, we can update the model parameters at no extra cost \citep{Aghamiry_2019_IWR}. This step is not mandatory for the algorithm. \\

\textbf{2) Estimating model parameters in $\mathbbm{2}$:} Because of the bi-linearity of the wave-equation \citep{Aghamiry_2019_IWR}, the medium parameters in $\mathbbm{2}$ are updated as 
\begin{equation} \label{m_2sub3}
\underset{\bold{m}_2\in \mathcal{M}}{\minimize} ~~~
\mathcal{R}_2(\bold{m}_2)+\sum_{i=1}^{n_s}\|\bold{L}(\bold{U}^{k+1}_{2,i})\bold{m}_2-\bold{y}_i^{k+1}\|_2^2,
\end{equation} 
where $\bold{U}^{k+1}_{2,i}$ refers to the wavefield in $\mathbbm{2}$ related to $i^{th}$ source,  $\bold{L}(\bold{U}^{k+1}_{2,i})=\omega^2 \text{Diag}(\bold{U}^{k+1}_{2,i})$, $\bold{y}_i^{k+1}=\bold{B}_{i}+\bold{\hat{B}}_{i}^k-\bold{A}_1(\bold{m}^{0}_{1})\bold{U}^{0}_{1,i}-{\Delta}_2\bold{U}^{k+1}_{2,i}$ and $\bold{B}_{i}$ and $\bold{\hat{B}}_{i}^k$ refer to the $i^{th}$ source and Lagrange multiplier, respectively 
and $\Delta_2$ is the Laplacian in $\mathbbm{2}$.  Eq. \eqref{m_2sub3} can be solved using splitting methods in the presence of non-smooth regularization terms \citep{Aghamiry_2019_CRO}.\\

\textbf{3) Estimating wavefields in $\mathbbm{2}$:} We update the wavefields in $\mathbbm{2}$ by finding the minimizer of Eq. \eqref{u2}, which has a closed form solution
\begin{align}
 \label{u_2sub1}
\bold{U}_2=&[\bold{A}_2(\bold{m}^{k+1}_2)^T\bold{A}_2(\bold{m}^{k+1}_2)]^{-1}[\bold{A}_2(\bold{m}^{k+1}_2)^T (\bold{B}+\bold{\hat{B}}^k-\bold{A}_1(\bold{m}^0_1)\bold{U}^{0}_1)].
\end{align}
The term $\bold{B}+\bold{\hat{B}}^k-\bold{A}_1(\bold{m}^0_1)\bold{U}^{0}_1$ is called redatumed data by \citet{Masson_2017_BTL}, because the sources and receivers virtually are moved onto the boundaries of $\mathbbm{2}$. The important point about LWI is that $\bold{U}^{0}_1$ in Eq. \eqref{u_2sub1} is reconstructed using the DA wavefield in Eq. \eqref{DA0}, which considers the data $\bold{D}$ for extracting the wavefield. This helps LWI to minimize the interaction of the wavefield between $\mathbbm{1}$ and $\mathbbm{2}$, when $\mathbbm{2}$ changes. The inclusion test (Fig. \ref{fig:FIG2}) is designed to show this effect.\\    

\textbf{4) Updating the dual variables:} At the final step, like IR-WRI, we update dual variables by the running sum of the constraint violations (source and data residuals), Eqs. \eqref{Bk}-\eqref{Dk}. \\
%
%
%
%
%
For a frequency batch with $n_b$ iterations, the proposed method requires $n_s$ PDE solutions of size $n \times n$ for step 2 and $n_b*n_s$ PDE solutions of size $n_2 \times n_2$ for step 5, while it is $n_b*n_s$ PDE solutions of size $n \times n$ in IR-WRI. The computational burden of other steps are negligible in LWI and IR-WRI. 
%
%
%
%
%
%
%
\section{Numerical results}
\begin{figure}
\center
\includegraphics[width=1.0\columnwidth]{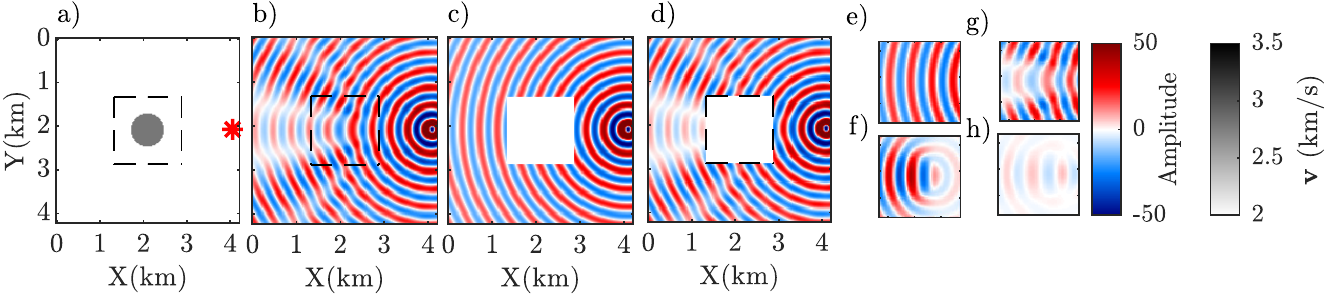}
\caption{(a) True model. Dashed lines indicate the boundary of $\mathbbm{2}$. The red asterisk indicates the location of the physical source. (b) 5 Hz wavefield calculated in (a). (c-d) Extracted wavefield of $\mathbbm{1}$ for homogeneous velocity model using (c) wave-equation, (d) Eq. \eqref{DA0}. (e) The extracted wavefield in $\mathbbm{2}$ with Eq. \eqref{u_2sub1} when wavefield (b) is used for $\mathbbm{1}$. (g) The difference between true wavefield (b) and (e). (g-h) Same as (e-f), but with DA wavefield shown in (d).}
\label{fig:FIG2}
\end{figure}

\textbf{Inclusion test:} 
When the model parameters change in $\mathbbm{2}$, the interaction of the wavefield between $\mathbbm{2}$ and $\mathbbm{1}$ are missed if we don't update the wavefield in $\mathbbm{1}$. This simple test is devoted to showing the effects of the DA wavefield $\bold{U}^{0}_1$ of $\mathbbm{1}$ in Eq. \eqref{u_2sub1} instead of using the wavefield extracted by solving the wave equation in the full domain.
The test involves a 4.2 km $\times$ 4.2 km homogeneous velocity model with $\bold{v}=2$~km/s and a circular inclusion with $\bold{v}=2.8$~km/s in the center depicted in Fig. \ref{fig:FIG2}a, where the boundary between $\mathbbm{1}$ and $\mathbbm{2}$ is shown with dashed lines. The acquisition consists of 120 receivers along a ring with a radius of $1.9$~km and a source shown by a red asterisk in Fig \ref{fig:FIG2}a. The 5 Hz wavefield related to this model is shown in Fig. \ref{fig:FIG2}b. 
Suppose we want to extract the wavefield in $\mathbbm{2}$ when the wavefield in $\mathbbm{1}$ (and the values on the boundaries of $\mathbbm{2}$) is calculated in the homogeneous medium with $\bold{v}=2$~km/s, without the inclusion. The extracted wavefields in $\mathbbm{1}$ by solving wave-equation and \eqref{DA0}, are shown in Figs. \ref{fig:FIG2}c-\ref{fig:FIG2}d, respectively. 
Because the DA wavefield involves the data term when solving the wave equation, it is a better approximation of the true wavefield in $\mathbbm{1}$ (compare wavefields in $\mathbbm{1}$ in Figs. \ref{fig:FIG2}b and \ref{fig:FIG2}d). 
Finally, the extracted wavefield in $\mathbbm{2}$ using \eqref{u_2sub1} with Figs. \ref{fig:FIG2}c-\ref{fig:FIG2}d are shown in Figs. \ref{fig:FIG2}e and \ref{fig:FIG2}g, respectively. The difference between these wavefields and the true wavefield in $\mathbbm{2}$ are shown in Figs. \ref{fig:FIG2}f and \ref{fig:FIG2}h, respectively. The extracted wavefield in $\mathbbm{2}$ using DA wavefield, Figs. \ref{fig:FIG2}g-\ref{fig:FIG2}h, can take into account the interactions between $\mathbbm{1}$ and $\mathbbm{2}$ and hence gives a better approximation of the wavefield in $\mathbbm{2}$. \\
%
%
%
%
%
%
\begin{figure}[t!]
\center
\includegraphics[width=0.6\columnwidth]{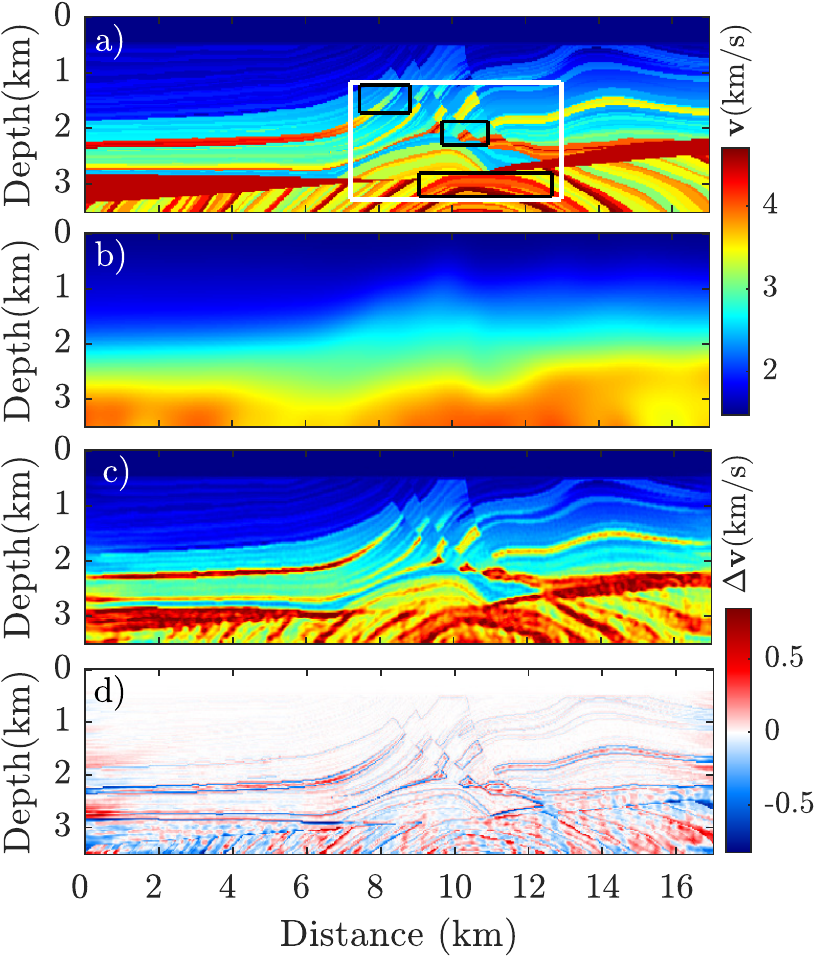}
\vspace{-0.25cm}
\caption{(a) True baseline model. The black lines indicate the boundaries of different parts of $\mathbbm{2}$. The white rectangle is the portion that is shown in Fig. \ref{fig:FIG4}. (b) The initial model for the inversion of baseline data. (c) Estimated IR-WRI velocity model applied on baseline data. (e) The difference between (a) and (c).}
\label{fig:FIG3}
\end{figure}
\textbf{4D FWI example:} 
Here the goal is to rapidly estimate the local changes that happen because of injected fluids or gas in the subsurface between a baseline and monitor data set \citep{Asnaashari_2012_TLI}.
In the past decade, different methods have been proposed for this problem. Among them, the most straightforward one is to apply conventional FWI on the monitor data when the recovered model from the baseline data is used as the initial model and then subtract the estimated models, which is expensive. To bypass this issue, one can use the velocity extracted from the baseline data as the initial model for applying localized FWI on the monitor data, which can decrease the computational burden significantly compared to the former case \citep{Malcolm_2016_R4F}.\\  
For this test, we use the Marmousi II benchmark (Fig. \ref{fig:FIG3}a) as the baseline model. The fixed-spread acquisition contains 57 point sources spaced 300~m apart and the source signature is a 10 Hz Ricker. A line of receivers spaced 50~m apart at the surface is used to record data. We apply the inversion in the 3~Hz - 13~Hz frequency band with a frequency interval of 1~Hz. We perform successive mono-frequency inversions proceeding from the low frequencies to the higher ones. We perform three passes through the frequencies, using the final model of one path as the initial model of the next one. The starting and final frequencies of the three paths are [3, 5], [3, 9], [5, 13]~Hz, respectively.
The initial model and the final estimated model using IR-WRI without regularization are shown in Figs. \ref{fig:FIG3}b-\ref{fig:FIG3}c, respectively. The difference between the true model and the final model of IR-WRI is shown in Fig. \ref{fig:FIG3}d. 
%
%
%
\begin{figure}[t!]
\center
\includegraphics[width=0.6\columnwidth]{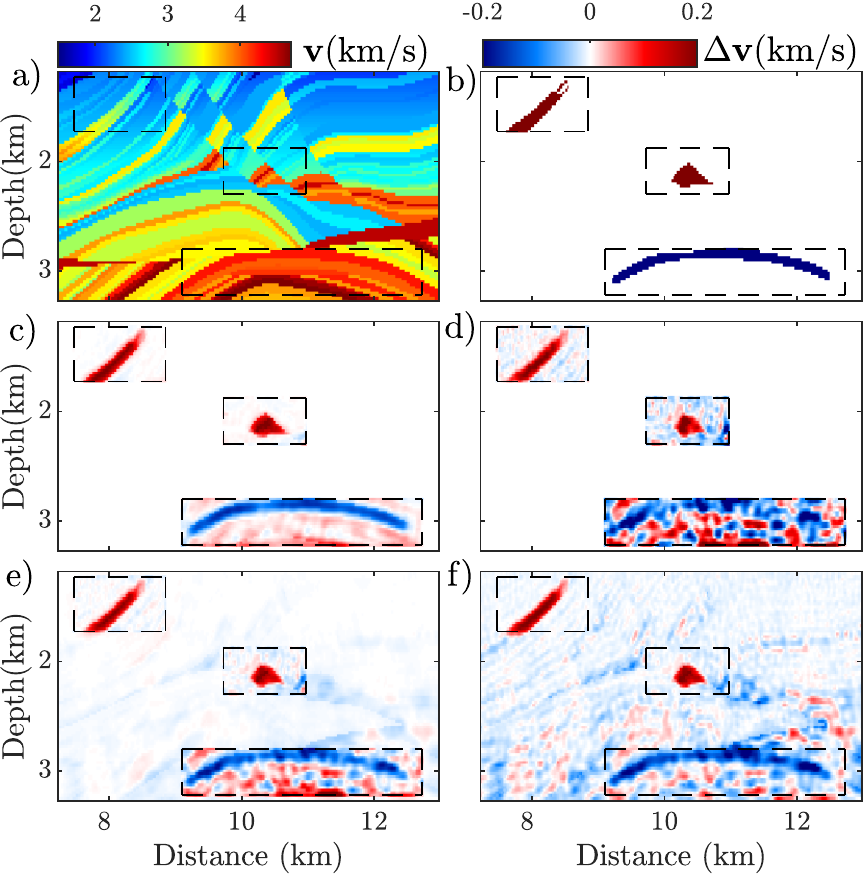}
\caption{(a) True monitor model. (b) The difference between the monitor and baseline model (Fig. \ref{fig:FIG3}a). (c) LWI with Fig. \ref{fig:FIG3}a as the initial model. (d) LWI with Fig. \ref{fig:FIG3}c as the initial model. (e) Same as (d), but with updating the background model. (f) IR-WRI with Fig. \ref{fig:FIG3}c as the initial model.}
\label{fig:FIG4}
\end{figure}
We create the monitor model by adding a small perturbation (Fig. \ref{fig:FIG4}a) to the true model. Also, the difference between monitor and baseline model is shown in Fig. \ref{fig:FIG4}b. Only the part of the Marmousi model, which is indicated by white rectangular in Fig. \ref{fig:FIG3}a is shown in Fig. \ref{fig:FIG4}. In this test $\mathbbm{2}$ contains three different parts as indicated by black rectangles in Figs. \ref{fig:FIG3}a and \ref{fig:FIG4}.  
We use three frequencies for the monitor data inversion, [5, 10, 15]Hz, and a successive mono-frequency inversion is applied with a stopping criteria of five iterations for each frequency. In addition, we perform two paths through the frequencies, which means that the total number of iterations is 30 for all the inversion results in Fig. \ref{fig:FIG4}. 
The extracted LWI model when the true baseline model (Fig. \ref{fig:FIG3}a) is used for the initial model of LWI is shown in Fig. \ref{fig:FIG4}c. We continue by applying LWI when the extracted model from IR-WRI of baseline data (Fig. \ref{fig:FIG3}c) is used as the initial model. These LWI results demonstrate the impact of the initial models used for the inversion with obvious differences in the largest, deepest rectangle. We repeat the LWI test in Fig. \ref{fig:FIG4}d, but we allow for updating the velocity model in $\mathbbm{1}$ one time without extra cost as discussed above (Fig. \ref{fig:FIG4}e). Comparing Figs \ref{fig:FIG4}d-\ref{fig:FIG4}e show significant improvements in the reconstructed model. Also, we apply IR-WRI (inversion in the entire of the medium) with the same setup and configuration of the LWI test on the monitor data (Fig. \ref{fig:FIG4}f). We observe that LWI and IR-WRI in $\mathbbm{2}$ (black rectangles in Figs. \ref{fig:FIG4}e-\ref{fig:FIG4}f) agree well.\\  
%
%
%
%
%
%
\section{Conclusion}
We implemented a localized version of the recently proposed IR-WRI, called LWI, to reduce computational cost for target-oriented applications.   
The proposed method does a block decomposition on the original FWI problem and decomposes the computational domain to a zone of interest $\mathbbm{2}$ and the rest ($\mathbbm{1}$). The new optimization problem is solved using an adaptation of multi-block ADMM when the subproblems related to $\mathbbm{1}$ are solved only once, but the others are solved normally. The proposed method is fast and straightforward, but it only approximately solves the subproblems related to $\mathbbm{2}$ when its velocities change. We show that this error is small because of the ability of the DA wavefield to approximate the true wavefield in $\mathbbm{1}$. Also, we show that the algorithm can update the background medium without additional computational cost when the wavefield of $\mathbbm{1}$ is updated. This update gives a better reconstruction in the targeted area when the velocity in $\mathbbm{1}$ is inaccurate. Numerical tests show promising results with this new algorithm.
%
%
%

\textbf{Acknowledgments:} 
H. Aghamiry would like to thank S. Beller (Geoazur) for fruitful discussions around this topic. A. Malcolm would like to thank the NSERC discovery grant program. This study was funded by the WIND consortium (\hyperlink{https://www.geoazur.fr/WIND}{www.geoazur.fr/WIND}). This study was granted access to the HPC resources of SIGAMM and the HPC resources of GENCI under the allocation 0596.
%
%
\bibliographystyle{seg}
\newcommand{\SortNoop}[1]{}

\end{document}